\documentclass{amsart}

\usepackage{amssymb}
\usepackage[all]{xy}
\usepackage{hyperref}

\usepackage{enumitem}   

\setlist[enumerate]{itemsep=.2em,topsep=.2em,leftmargin=1.25em,itemindent=2.0em}

\newtheorem{thm}{Theorem}

\newtheorem{lem}[thm]{Lemma}
\newtheorem{cor}[thm]{Corollary}

\newtheorem{prop}[thm]{Proposition}

\theoremstyle{definition}
\newtheorem{defn}[thm]{Definition}

\newtheorem{say}[thm]{}
\newtheorem{exmp}[thm]{Example}
\newtheorem{exmps}[thm]{Examples}

\newtheorem{rem}[thm]{Remark}
\newtheorem{rems}[thm]{Remarks} 

\newtheorem*{ack}{Acknowledgments}      
\newtheorem{notation}[thm]{Notation}   
  
\newtheorem{defn-thm}[thm]{Definition--Theorem}  
\newtheorem{defn-lem}[thm]{Definition--Lemma}  

\newtheorem{prevress}[thm]{Previous results}

\theoremstyle{remark}

\renewcommand{\c}[0]{{\mathbb C}}  

\renewcommand{\o}[0]{{\mathcal O}} 
\newcommand{\z}[0]{{\mathbb Z}}

\renewcommand{\a}[0]{{\mathbb A}}

\newcommand{\p}[0]{{\mathbb P}}

\newcommand{\q}[0]{{\mathbb Q}}
\newcommand{\map}[0]{\dasharrow}
\newcommand{\qtq}[1]{\quad\mbox{#1}\quad}
\newcommand{\spec}[0]{\operatorname{Spec}}

\newcommand{\pics}[0]{\operatorname{\mathbf{Pic}}}

\newcommand{\pico}[0]{\operatorname{\mathbf{Pic}}^{\circ}}

\newcommand{\ns}[0]{\operatorname{NS}}

\newcommand{\totaldiscrep}[0]{\operatorname{totaldiscrep}}

\newcommand{\supp}[0]{\operatorname{Supp}}    
\newcommand{\red}[0]{\operatorname{red}}    
\newcommand{\codim}[0]{\operatorname{codim}}    
    
\newcommand{\proj}[0]{\operatorname{Proj}}

\newcommand{\cent}[0]{\operatorname{center}}

\newcommand{\diff}[0]{\operatorname{Diff}}

\newcommand{\simq}[0]{\sim_{\q}}

\newcommand{\depth}[0]{\operatorname{depth}}

\def\into{\DOTSB\lhook\joinrel\to}

\def\loccoh#1.#2.#3.#4.{H^{#1}_{#2}(#3,#4)}

\DeclareMathAlphabet{\mathchanc}{OT1}{pzc}%
                                {m}{it}

\newcommand{\gm}[0]{{\mathbb G}_m}
\newcommand{\ga}[0]{{\mathbb G}_a}

\usepackage[all]{xy}\xyoption{dvips}

\newcommand{\auts}[0]{\operatorname{\mathbf{Aut}}}

\newcommand{\dom}[0]{\operatorname{domi}}

\newcommand{\ner}[0]{\operatorname{N\acute{e}r}}

\newcommand{\nero}[0]{\operatorname{N\acute{e}r}^{\circ}}

\begin{document}
\bibliographystyle{amsalpha}


 \title[N\'eron models and   minimal models]{N\'eron models,  minimal models, and \\ birational group actions}
 \author{J\'anos Koll\'ar}

 \begin{abstract}
   Let $A_K$ be an Abelian variety over the quotient field of a Dedekind domain $R$.
   We show that the identity component of the N\'eron model of $A_K$ acts regularly on any minimal model of $A_K$ over $\spec R$, and discuss when the N\'eron model is an open subset of a minimal model.

   The main technical result is the rigidity of small modifications, leading to a  regularity criterion for  birational group actions.

       \end{abstract}

  \maketitle

Let $S$ be a 1-dimensional, irreducible, regular scheme with function field $K$,
and $A_K$ an Abelian variety over $K$.

Finding a model of $A_K$ over $S$ attracted considerable attention.
   Kodaira described minimal,  elliptic  surfaces  over $\c$
\cite{MR132556, MR184257}; 
        this was extended to arbitrary Dedekind domains in
       \cite[Chap.III]{MR179172}.
            For higher dimensional Abelian varieties,   the optimal  smooth group scheme  $\ner(A_K)$---called the {\it N\'eron model}---is constructed in  \cite[Chap.II]{MR179172};
       see  Section~\ref{ner.sec} for a summary.

       Minimal model theory---also known as Mori's program---aims to  provide a (non-unique) projective, minimal  model  of any variety or scheme; see  (\ref{main.thm.rems}.4)  and (\ref{min.mod.say}). For Abelian generic fibers,
         minimal  models give higher dimensional analogs of  the  minimal,  elliptic surfaces of Kodaira and N\'eron.
       
         The aim of this note is to study
actions of group schemes on minimal models, and then derive 
 connections between the
 N\'eron model and the minimal models.
 The main result is Theorem~\ref{main.min.ner.thm.gen}, which
 gives the following for the N\'eron model.

  \begin{thm} \label{main.min.ner.thm}
    Let $S$ be a 1-dimensional, irreducible, regular scheme with function field $K$, and $A_K$ an Abelian variety over $K$. Assume that  $A_K$ has a  projective, minimal  model  $\pi: A_S\to S$.
    Let $\nero(A_K)\subset \ner(A_K)$ be the identity component of the N\'eron model (\ref{nero.defn}).
    
     Then   the $A_K$-action on itself extends to a regular $\nero(A_K)$-action on $A_S$.
        \end{thm}

  Let $p:Y_S\to S$ be a morphism with generic fiber $A_K$, and
  $Y_S^{\rm sm}\subset Y_S$ the $p$-smooth locus. By definition, there is a birational morphism   $Y_S^{\rm sm}\to \ner(A_K)$, which is an isomorphism if $Y_S$ is a minimal, elliptic surface  by \cite[Chap.III]{MR179172}.
  In higher dimensions the induced
  $A_S^{\rm sm}\to \ner(A_K)$
is an isomorphism in some cases.

  \begin{cor} \label{main.min.ner.thm.cor}
    Let $\pi: A_S\to S$ be a projective, minimal model. Then
    the natural morphism  $A_S^{\rm sm}\to \ner(A_K)$ is an
  \begin{enumerate}
  \item    open embedding, whose image is $\ner(A_K)\setminus(\mbox{finitely many fibers})$, and an
\item       isomorphism   in  either of the following cases:
 \begin{enumerate}
 \item $A_S$ is regular, or
 \item $A_S$ is factorial and the residue fields of $S$ are  perfect.
   \end{enumerate}
  \end{enumerate}
  \end{cor}

  \begin{prevress} There is a vast literature on
    degenerations of Abelian varieties, focusing on the case of
    semiabelian reduction, which can be always be achieved
    after  a suitable finite base change  $S'\to S$; see  \cite{MR0352106, MR1083353, MR1646554,  MR1707764, ale-abvar, mitsui2024relativecompactificationssemiabelianneron, nakamura2024relativecompactificationsemiabelianneron} and the references there. However, the situation seems quite different without a base change.

    The  case (\ref{main.min.ner.thm.cor}.2.a)  is 
    treated in \cite{MR3770167, MR3826826, bkv}; see (\ref{k.m.n.m.say}) for details and comments.

    The group of birational automorphisms of a projective minimal model is studied in \cite[3.3]{hanamura}.
     A  detailed description of  N\'eron models of  Lagrangian fibrations over higher dimensional bases is given in \cite{kim2024neronmodelhigherdimensionallagrangian}. In particular, \cite[Thm.1.3]{kim2024neronmodelhigherdimensionallagrangian}
proves that the identity component of the N\'eron model acts regularly on the total space of a Lagrangian fibration.
  \end{prevress}

  In Theorem~\ref{main.min.ner.thm}
and Corollary~\ref{main.min.ner.thm.cor} we assume that a minimal model exists, but after that the proofs do not use  any minimal model techniques.
It is sufficient to know that something called a minimal model with generic fiber $A_K$ exists, it is normal, and 
(\ref{min.mod.say}.4) holds. 

One can get further results  relying on the construction of  additional minimal models involving  boundary divisors, and using more information about the occurring  singularities. Currently these are known in  characteristic 0 (\ref{main.thm.rems}.4).
The main result of this type is the extension of
Theorems~\ref{main.min.ner.thm} and~\ref{main.min.ner.thm.gen} to
log canonical pairs, given  in (\ref{logmin.lc.thm}) and (\ref{gp.sch.acts.mm.prop}).

        Following \cite{MR715605}, a {\it good completion of a N\'eron model}
        is an open embedding $\ner(A_K)\into Y$ such that $Y$ is normal, projective over $S$, and the induced rational $\ner(A_K)$-action  on $Y$ is regular;
        {\it equivariant completion} may be a better name.
        Their existence is stated as
        conjecture   GCNM   on \cite[p.318]{MR715605}, and attributed to Mumford. We prove  it  in Section~\ref{eq.comp.sec} in characteristic 0.

\begin{cor} \label{main.min.ner.thm.cor.2}
  If $S$ is   over a field of characteristic 0, then
   \begin{enumerate}
   \item there is a minimal model $A^*_S$ on which $\ner(A_K)$ acts regularly, and 
   \item $\ner(A_K)\to S$ has an equivariant completion.
      \end{enumerate}
\end{cor}

The papers 
\cite{mitsui2024relativecompactificationssemiabelianneron, nakamura2024relativecompactificationsemiabelianneron}
construct an equivariant completion of the N\'eron model 
for  the semiabelian reduction cases, but  do not check that its singularities are terminal. Using  minimal model  theory, we get the following.

        \begin{cor} \label{main.min.ner.thm.cor.3}
          Assume that all  fibers of  $\ner(A_K)\to S$ are semiabelian, and $S$ is   over a field of characteristic 0. Then
          \begin{enumerate}
 \item all fibers of $\pi: A_S\to S$ are reduced, semi-log-canonical,  
 \item         $\ner(A_K)$ is isomorphic to the $\pi$-smooth locus,
   \item for every  $s\in S$, $\ner(A_K)_s$ acts with a dense, open orbit on $A_s$, and
\item   $A_s$ is  {\em semiabelic}  in the terminology of \cite[1.1.5]{ale-abvar}.
   \end{enumerate}
        \end{cor}

        \begin{rems}\label{main.thm.rems} {\ }
            
          (\ref{main.thm.rems}.1)  Minimal models are not unique, and our methods establish properties that hold for {\em all} of them.
          Fixing  a suitable polarization on $A_K$ may lead to  a unique model $A_S$. The semiabelian cases are worked out in 
          \cite{MR0352106, MR1083353, MR1646554,  MR1707764, ale-abvar, mitsui2024relativecompactificationssemiabelianneron, nakamura2024relativecompactificationsemiabelianneron}.

(\ref{main.thm.rems}.2)
 Usually there is no regular $\ner(A_K)$-action on a given minimal model; see
 (\ref{many.minmods.exmp}) for some examples.

          (\ref{main.thm.rems}.3)  Let $\sum_{i\in I} m_iE_i$ be a fiber of $\pi$.
          If $\dim A_K=1$, then either   $m_i=1$ for some $i$, or
          $\gcd(m_i: i\in I)>1$. I do not know whether this holds in higher dimensions.

          (\ref{main.thm.rems}.4)
          Roughly speaking, minimal models are known to exist
        in characteristic 0,  and are conjectured to exist in general.
       In characteristic 0 the main references are
        \cite{MR2359343, bchm, lyu-mur}, while 
  \cite{many-p} settles the
         case   $\dim A_K\leq 2$, $S$  excellent, and all residue characteristics are $>5$.  For Abelian generic fibers, see also \cite{fujino-ssmmp, lai-2009} for key steps of the proof.

(\ref{main.thm.rems}.5)
         By Nagata's theorem, every separated $S$-scheme $X$ of finite type  is isomorphic to an  open subscheme of a proper $S$-scheme  $\bar X$; see  \cite{MR0158892}.  Over a field, equivariant completions are established in
         \cite{MR337963,  MR4504264}.
           The case of  smooth, affine group schemes  is considered in   \cite{MR387294}, but it is not clear  whether all the assumptions there are necessary. 
         Very little seems to be known about equivariant completions of actions of non-affine  group schemes.

         (\ref{main.thm.rems}.6) A stronger form of the existence of
         equivariant completions is conjectured in
         \cite[12.4]{nakamura2024relativecompactificationsemiabelianneron}.

         (\ref{main.thm.rems}.7)  Y.-J.~Kim pointed out that $\nero(A_K)$ is isomorphic to the main component of  $\auts_S(A_S)$ using 
         \cite[2.42]{kim2024neronmodelhigherdimensionallagrangian}.

 \end{rems}

        \begin{ack}  I thank V.~Alexeev,  Y.-J.~Kim, R.~Laza, S.~Molcho and I.~Nakamura for  comments and references.
  Partial  financial support    was provided  by  the NSF (grant number
DMS-1901855)  and by the Simons Foundation   (grant number SFI-MPS-MOV-00006719-02).
\end{ack}

\section{Main technical theorems}\label{main.thm.sec}

Theorem~\ref{main.min.ner.thm} is a special case of the following regularity criterion, whose proof is given in (\ref{main.min.ner.thm.gen.pf}). 
It can be viewed as a generalization of
 \cite[3.3]{hanamura} and of \cite[Thm.1.3]{kim2024neronmodelhigherdimensionallagrangian}.
      
  \begin{thm} \label{main.min.ner.thm.gen}
    Let $S$ be a  normal scheme with function field $K$, 
     $G_S\to S$ a smooth group scheme with connected fibers, and 
     $\pi: X_S\to S$  a dominant, projective morphism.
    Assume that $X_S$ has terminal singularities and $K_{X_S/S}$ is $\pi$-nef.
    
    Then any rational  $G_K$-action
    $\mu_K: G_K\times_K X_K\map X_K$ on the generic fiber   $X_K$
    extends to a  regular  $G_S$-action
    $\mu_S: G_S\times_S X_S\to X_S$.
  \end{thm}

  Note that the statement tacitly assumes that $K_{X_S/S}$ exists.
  This holds if $S$ is excellent and has a dualizing complex.

  Theorem~\ref{main.min.ner.thm.gen} is derived from a more general regularity criterion  on
  rational group actions   (\ref{gp.sch.acts.prop}). The full scheme theoretic version is a little technical, the following variant for varieties, to be proved in (\ref{sm.mod.gr.isom.cor.say}), gives the right idea.

\begin{thm}\label{sm.mod.gr.isom.cor}
  Let $X$ be a normal, projective variety  over an algebraically closed field, and $G$ a connected, smooth algebraic group. Let $\mu:G\times X\map X$ be a rational action. Assume that, for general $g\in G$, the map $\mu_g:X\map X$ is an isomorphism in codimension 1.
  Then  $\mu:G\times X\to X$ is a regular action.
\end{thm}

Note that we do not assume that there is an open subset $U\subset X$ where all the $\mu_g$ are defined.

  \begin{say}[Small modifications]\label{small.mod.int.say}
    Let $X_i$ be normal, proper schemes over $S$. A birational map
    $\phi:X_1\map X_2$ is a {\it small modification,}
    if there are closed subsets 
   $Z_i\subset X_i$ of codimension $\geq 2$ such that   $\phi$ restricts to an isomorphism
    $\phi:(X_1\setminus Z_1)\cong (X_2\setminus Z_2)$.

    The best known examples are flips and flops, much studied when  the $X_i$ have log terminal singularities \cite{km-book}. The case of arbitrary normal varieties attracted much less attention  \cite[Sec.2]{k-etc}. 

    The key technical result behind Theorem~\ref{sm.mod.gr.isom.cor} is the rigidity of small modifications, see (\ref{small.mod.prop}) and (\ref{small.mod.prop.gen}) for the precise statements.

    Very little seems to be known about the  existence and structure of
    small modifications, beyond the log terminal case.
    As the next example shows, there are some unexpected  flops of 3-folds with log canonical singularities.
    \end{say}

  \begin{exmp}[Flopping a fiber of a product] \label{small.mod.prop.2.exmp}
    We construct an example of a pointed surface $v\in S$ with simple elliptic (hence log canonical) singularity at $v$, and an elliptic curve $E$ such that the fiber
    $E\times \{v\}\subset E\times S$ has infinitely many different flops.

  Let $E\subset \p^2$ be a smooth cubic,
  $S\subset \a^3$ the cone over $E$ with vertex $v$ and  projection $\pi:S\map E$.
  For $m\in\z$ consider the map
  $$
  \phi_m: E\times S\map E\times S\qtq{given by} (p, s)\mapsto \bigl(p+m\pi(s), s\bigr).
  $$
  Then $\phi_m$ is an isomorphism over $S\setminus \{v\}$, but it is not defined along $E\times \{v\}$ for $m\neq 0$. 
  Note that $E\times S$ has log canonical singularities, and the $\phi_m$ give infinitely many different flops of the curve $E\times \{v\}$.
  \end{exmp}

        \section{N\'eron models}\label{ner.sec}

        We summarize the properties of N\'eron models that we use;
        see \cite{MR861977, blr} for complete treatments, and
        \cite{MR3695053} for futher results.

\begin{notation}\label{dvr.notation}
  As explained in \cite[p.13]{blr}, the theory of N\'eron models is local on the base. Thus, from now on, 
          $R$ denotes a DVR with quotient field $K$, residue  field $k$ and $T:=\spec R$.
  For a morphism  $\pi:X_T\to T$, we use $X_K$ (resp.\ $X_k$) to denote the
  generic (resp.\ closed or central) fiber.
         \end{notation}

       \begin{defn} \label{ner.mod.defn}
  A smooth morphism $X_T\to T$ satisfies the
  {\it  N\'eron mapping property} iff the following holds.
 \begin{enumerate}
    \item
  Let $Y_T\to T$ be any smooth morphism. Then every morphism  $\phi_K: Y_K\to X_K$ extends to a morphism  $\phi_T: Y_T\to X_T$.
\end{enumerate}
 A   morphism $X_T\to T$  is a {\it N\'eron model} of its generic fiber $X_K$, if it is smooth, of finite type, and satisfies the N\'eron mapping property. We use  $\ner(X_K)$   to denote  the  N\'eron model of $X_K$.

 It is clear that $X_K$ uniquely determines its N\'eron model $\ner(X_K)$.
 N\'eron models commute with 
 \'etale base changes  $T'\to T$  \cite[p.13]{blr}.
 Thus many questions about N\'eron models can be treaded over
 the strict Henselisation of $T$.

 However,  ramification usually changes
 the  N\'eron model drastically.
       \end{defn}

       \cite{MR146035, MR179172} construct N\'eron models  when $X_K$ is a
            {\it torsor} (or {\it principal homogeneous space})
            under an Abelian variety over $K$.
Note that 
\cite{MR861977} and \cite{blr}
discuss mainly (\ref{ner.mod.say}.1). 

\begin{thm}[N\'eron models]\label{ner.mod.say}
  Let $A_K$ be an  Abelian variety over $K$, and $X_K$  an $A_K$-torsor.
  Then the following hold.
   \begin{enumerate}
    \item $\ner(A_K)$ exists, and it  is a
          smooth, finite type  group scheme over $T$.
        \item  $\ner(X_K)$  exists.
        \item  $\ner(X_K)$ is either $X_K$ or a $\ner(A_K)$-torsor.
            \qed
  \end{enumerate}
\end{thm}

\begin{defn}\label{nero.defn}
   The union of   $A_K$ and the identity component of $\ner(A_K)_k$ 
  is  the {\it identity component} of  $\ner(A_K)$,  denoted by $\nero(A_K)$.
  Note that $\nero(A_K)$ is open and dense in $\ner(A_K)$, and
  $\nero(A_K)\to T$  is smooth group scheme with connected fibers.

  The central fiber $\ner(X_K)_k$ is usually disconnected.
  Its {\it group of connected components}   is a  finite abelian group.
    \end{defn}

\begin{say}[Pseudo-Abelian varieties]\label{pseu-ab.say}
  When $K$ is not perfact, it is natural to also consider
  pseudo-Abelian varieties $P_K$ as defined in 
  \cite{MR3185350}. These also have N\'eron models $\ner(P_K)$ by
  \cite[Sec.10.2]{blr}; see also 
  \cite{MR4784847} for further results.
 Theorem~\ref{main.min.ner.thm.gen} also applies to
  actions of  $\nero(P_K)$. 
\end{say}

 \section{Minimal models and N\'eron models}\label{minmod.sec}

 An introduction to minimal models is in \cite{km-book}, and
 a summary of later results is in \cite[Secs.11.1--2]{k-modbook}.
 For the purposes of  Theorem~\ref{main.min.ner.thm} we need only the following special setting.

 \begin{say}[Minimal models of Abelian torsors]\label{min.mod.say}
  
Let $S$ be a 1-dimensional, irreducible, regular scheme with function field $K$.
  Let $X_K$ be a proper $K$-variety with  terminal singularities, whose canonical class  is numerically trivial.
   
        A  proper morphism  $\pi:X_S\to S$ is a {\it minimal model}  of $X_K$ if 
         \begin{enumerate}
       \item $X_S$ has terminal singularities \cite[2.34]{km-book},
       \item  $K_{X_S}$ is numerically $\pi$-trivial, and
         \item the generic fiber is isomorphic to $X_K$.
         \end{enumerate}
{\it Numerically $\pi$-trivial} means that $(K_{X_S}\cdot C)= 0$ for every curve contracted by $\pi$.

           Being terminal  is preserved by \'etale base change \cite[2.15]{kk-singbook}, so being a minimal model is also preserved by \'etale base change.
       Note that minimal model theory frequently works with     {\it $\q$-factorial} schemes,  that is, some multiple of any Weil divisor is required to be Cartier.  However, being $\q$-factorial is not an \'etale-local property, so, as in \cite{MR4342251},  we  work without it.

          \medskip

Minimal models of $X_K$ need not be unique, but the following 
       is easy to see directly; see \cite[3.52]{km-book}.

             \medskip

             {\it Claim \ref{min.mod.say}.4.}
             Let $B$ be a  normal scheme  and 
     $\pi_i: X_i\to B$   proper morphisms.
    Assume that  $X_i$ is irreducible, has terminal singularities, and $K_{X_i/B}$ is $\pi_i$-nef for $i=1,2$.
    Let   $\phi:X_1\map X_2$  be a birational map.

    Then there are   closed $Z_i\subset X_i$ of codimension $\geq 2$, such that
    $\phi$ induces an isomorphism $X_1\setminus Z_1\cong X_2\setminus Z_2$.
            \qed

           \medskip

            For the existence  of minimal models, the main references are 
\cite{MR2359343, bchm}, with  additional steps   in \cite{fujino-ssmmp, lai-2009, lyu-mur}.

       \medskip

           {\it Theorem \ref{min.mod.say}.5.}  Assume that $S$ is over a field of characteristic 0.
         Let $X_K$ be a projective $K$-variety with  terminal singularities, whose canonical class  is numerically trivial.
         Then $X_K$ has a $\q$-factorial, projective,   minimal model
         $\pi:X_S\to S$.
         
         In particular,  every
         Abelian torsor over $K$ has a minimal model.
         \qed
       \end{say}

    The following generalization of 
    Theorem~\ref{main.min.ner.thm}
    connects  minimal models and N\'eron models.
    We follow Notation~\ref{dvr.notation}, thus $T$ is the spectrum of
    a DVR  $R$ with quotient field $K$ and residue  field $k$.

\begin{thm} \label{G.circ.acts.thm}
  Let  $X_K$  be a torsor under an  Abelian variety $A_K$,
  and  $\pi:X_T\to T$  a minimal model  of $X_K$.
  Let $\nero(A_K)\subset \ner(A_K)$ be the identity component of the N\'eron model.  Then the $A_K$-action on $X_K$ extends to a regular $\nero(A_K)$-action on $X_T$.
\end{thm}

  Proof.           The $A_K$-action on $X_K$ extends to a rational $\nero(A_K)$-action on $X_T$ which we write as 
          $$
  \mu: \nero(A_K)\times_T X_T\map X_T.
  \eqno{(\ref{G.circ.acts.thm}.1)}
          $$
          Being an morphism can be checked after a  surjective, \'etale base change.          We may thus assume that $T$ is strictly Henselian.

Given any section  $\tau:T\to \nero(A_K)$, the induced map
$\mu_{\tau}: X_T\map X_T$ is a small modification by
(\ref{min.mod.say}.4). Thus $\mu$ is a morphism by (\ref{gp.sch.acts.prop}).  
\qed
\medskip

The $\nero(A_K)$ action gives strong restrictions on the smooth locus.

\begin{lem} \label{sm=orbit.lem}
  Let $U_T\to T$ be a smooth morphism with   a  regular  $\nero(A_K)$ action,
  whose generic fiber $U_K$ is an $A_K$-torsor.
  Then $u:U_T\to \ner(U_K)$ is an open immersion, and the complement of the image is a union of 
 $\nero(A_K)_k$-orbits.
\end{lem}

Proof.  By the universal property of  N\'eron models there is a birational
morphism
$u:U_T\to \ner(U_K)$, which is $\nero(A_K)$-equivariant.  The restriction of $u$ gives a $\nero(A_K)_k$-equivariant morphism of   $U_k$ to  $\ner(U_K)_k$, which is a union of $\nero(A_K)_k$-orbits.  Thus $u_k$ is  an isomorphism onto an union of connected components of  $\ner(X_K)_k$. \qed
\medskip

The following is a generalization of
Corollary~\ref{main.min.ner.thm.cor}.

\begin{cor} \label{G.circ.acts.thm.cor} Continuing with the notation of
  (\ref{G.circ.acts.thm}), let $X^{\rm sm}_T\subset X_T$ denote the 
  $\pi$-smooth locus, and $u:X^{\rm sm}_T\to \ner(X_K)$ the induced morphism.
  Then
  \begin{enumerate}
  \item either $u$ is an isomorphism,
  \item or $X^{\rm sm}_T\cong X_K$.
      \end{enumerate}
  Furthermore, $u$ is an isomorphism in either of the following cases.
 \begin{enumerate}\setcounter{enumi}{2}
 \item $X_T$ is regular.
   \item $X_T$ is factorial and $k$ is perfect.
      \end{enumerate}
\end{cor}

Proof.  
 We may  assume that $T$ is strictly Henselian.
If $X^{\rm sm}_k$ is empty, we have
(\ref{G.circ.acts.thm.cor}.2).
 Otherwise there is at least one irreducible component
 $E_1\subset X_k$, such that $u$ maps  $E_1\cap X^{\rm sm}_k$
 onto an irreducible component $F_1\subset \ner(X_K)_k$.

 Let  $F_2\subset \ner(X_K)_k$ be any other irreducible component.
 There is a section $\tau:T\to  \ner(A_K)$ that maps $F_1$ onto $F_2$.
 This gives a rational map  $\mu_{\tau}:X_T\map X_T$ that sends
 $E_1$ to another irreducible component $E_2\subset X_k$ by
 (\ref{min.mod.say}.4).
 Then $u$ maps  $E_2\cap X^{\rm sm}_k$
 onto $F_2$. Thus (\ref{G.circ.acts.thm.cor}.1) holds by  (\ref{sm=orbit.lem}).

 It remains to show that $X^{\rm sm}_k$ is not empty in cases  (\ref{G.circ.acts.thm.cor}.3--4).
Let $C\subset X_T$ be a section. 
Then  $(X_k\cdot C)=1$, and
$(\red X_k\cdot C)\geq 1$ if $\red X_k$ is a Cartier divisor.
Thus $X_k$ is reduced near $X_k\cap C$, hence
geometrically reduced if $k$ is perfect.

If $X_T$ is regular, then $(X_k\cdot C)=1$ implies that $X_k$ is regular
at the $k$-point $X_k\cap C$, hence smooth. These show (\ref{G.circ.acts.thm.cor}.3--4).
 \qed

\medskip
The next examples show that in general $\nero(A_K)$
is not isomorphic to an open subset of any minimal model, and that minimal models are not unique.

\begin{exmp}\label{non.red.fiber.exmp.1}
 Let $B$ be an Abelian variety and $\sigma$ an automorphism of order $m$ fixing the identity.  Set
 $
 X:=X(B,\sigma):=B\times \a^1/(\sigma,\epsilon),
 $
 where  $\epsilon$ is  a primitive $m$th root of unity.
The central fiber has multiplicity $m$, and the reduced central fiber is  the Kummer-type variety  $B/\langle \sigma\rangle$.
The singularities of $X$ are analytically quotient singularities. The Reid-Thai criterion determines whether they are terminal \cite{r-c3f}.
If they are, then $X\to \a^1$ is a minimal model which does not contain
the N\'eron model. Let us see some concrete examples.

\smallskip
(\ref{non.red.fiber.exmp.1}.1)  Let $\sigma:=-1$ and $n=\dim B$.
The singularities of $X$ are analytically of the form  $\a^{n+1}/\frac12(1,\cdots,1)$.  These are  terminal for $n\geq 2$, and  $K_X\sim 0$ iff $n$ is odd.
 The    central fiber of the N\'eron model  is  $\ga^n\times (\z/2)^{2n}$.

\smallskip
(\ref{non.red.fiber.exmp.1}.2)  Set $E=(x^3+y^3+z^3=0)$ with $\z/3$-action
$(x{:}y{:}z)\mapsto (\epsilon x{:}y{:}z)$. Set $B=E\times E$ and
let $\z/3$-act as $(\epsilon, \epsilon, \epsilon^2)$.
There are 9 fixed points, all are analytically of the form  $\a^{3}/\frac13(1,1,2)$.  These are  terminal.

\smallskip
(\ref{non.red.fiber.exmp.1}.3)  With  $B=E\times E$ as above, 
let $\z/3$-act as $(\epsilon, \epsilon, \epsilon)$.
There are 9 fixed points, all are analytically of the form  $\a^{3}/\frac13(1,1,1)$.  These are  canonical but not terminal. A minimal model is obtained by blowing them up. We have now 9 reduced divisors in the cetral fiber, which cointains the  N\'eron model.
The    central fiber of the N\'eron model  is  $\ga^2\times (\z/3)^{2}$.

\smallskip
        (\ref{non.red.fiber.exmp.1}.4)  Let $B$ be a simple abelian variety of dimension $\geq 4$. Then  $B/\langle\sigma\rangle$ has canonical singularities by \cite[Thm.9]{k-lars}, hence $X(B,\sigma)$ has terminal singularities, and it is a minimal model.
 \end{exmp}

  \begin{exmp}\label{many.minmods.exmp}
    Let $g_i(x,y,z)$ be 2 general cubics and
     $\pi:X\to \a^1$  the fiber product of the 2 elliptic surfaces
    $\bigl(xyz+tg_i(x,y,z)=0\bigr)\subset \p^2\times \a^1$.

    $K_X$ is trivial and the general fiber is a product of 2 elliptic curves. The central     fiber is a union of 9 copies $L_{ij}\cong\p^1\times \p^1$.
    $X$ has  9 singular points  of the form $(xy=uv)$ where 4 of the $L_{ij}$  meet.
    The central     fiber of the N\'eron model is  $\gm^2\times (\z/3)^2$.
    It acts regularly on $X$, transitively on both the $L_{ij}$ and the singular points.

    Each singular point has 2 small resolutions, giving $3^9$ different
    minimal models, among them $2^9$ are smooth. The projective ones  can be obtained by blowing up a divisor  $\sum_{ij} m_{ij}L_{ij}$. For example,
    $\sum_i L_{ii}$  gives  a projective,  smooth minimal model.

    There are 2 smooth minimal models on which the whole  N\'eron model  acts, but neither is projective.
  \end{exmp}

  \begin{exmp}\label{sing.minmods.exmp}  As a variant of
    (\ref{many.minmods.exmp}), work over the base field $\q$ and
    start with the elliptic surfaces 
    $$
    S_a:=\bigl(x^3+ay^3+a^2z^3-3axyz=tg(x,y,z)\bigr),
    $$
    whose central fiber  is a union of the 3 lines
    $\bigl(x+\alpha y+\alpha^2z=0\bigr)$ for $\alpha^3=a$.

    Now consider the fiber product of $S_a$ and $S_c$, say for $a=2, c=3$.
    The central fiber  is the union of 9 divisors, which are conjugate to each other. It has 9 singular points, and no projective, smooth, small  modifications.
    \end{exmp}

  \begin{say}[Kulikov models]\label{k.m.n.m.say}
    Working over $\c$, 
the papers 
 \cite{MR3770167, MR3826826, bkv} state
 that the $\pi$-smooth locus of a  {\it Kulikov model}   $\pi: \operatorname{Kul}(A_K)\to T$ of an Abelian variety $A_K$ is
 isomorphic to   its N\'eron model.

         Note however,  that
\cite{MR3826826, bkv} and   \cite{MR3770167}
          use non-equivalent definitions of a
          Kulikov model.

         In all 3  papers, a Kulikov model is a flat, proper morphism
         $\pi: Y\to T$,  where $Y$ is regular, and the general fiber
         has trivial canonical class  (in our case,  an Abelian variety).
         The additional assumptions are, however, quite different.

         In \cite[5.1.2]{MR3770167} the central fiber $Y_0$  is required to be a  normal crossing  divisor, and the log canonical class
         $K_{Y/T}(\red Y_0)$ is assumed to be numerically  $\pi$-trivial. 

         \cite[2.26]{bkv} refers to \cite[3.3]{MR3826826}.
         In the latter the central fiber of a Kulikov model is allowed to have arbitrary singularities,  and
         $K_{Y/T}$ is assumed to be numerically  $\pi$-trivial; see  \cite[3.2]{MR3826826}.

          The arguments in   \cite{MR3826826, bkv}   rely on   \cite{MR3770167},
          sliding over       the incompatibility of the definitions.

 Note also that \cite[3.3]{MR3826826} refers to 
 \cite[6.7]{MR3770167}. There does not seem to be  any entry numbered (6.7) in
 \cite{MR3770167}; neither in  the published version nor in any of the three ArXiV versions.
\end{say}

\section{Small modifications}

We work over a Noetherian base scheme. To avoid some possible problems with dimension, let us assume that it is universally catenary.

\begin{defn}[Families of birational maps] \label{bir.fam.defn}
A {\it birational map}  $\phi:X\map Y$ between (reduced, separated) schemes $X,Y$ is an  (equivalence class of) 
  isomorphisms between dense, open subsets  $ X\supset U_X\cong U_Y\subset  Y$.
  The largest such is called the {\it domain of isomorphy} of $\phi$, denoted by
  $\dom(\phi)$, and usually viewed as a subset of $X$.
  
Let $p_i:X_i\to S$  be  equidimenisonal, proper morphisms.
A birational map
$\phi: X_1\map X_2$  over $S$ is a
{\it family of birational maps} if
$\dom(\phi)$ contains  the generic points of each fiber. Thus the restrictions
$\phi_s: \red(X_1)_s\map \red(X_2)_s$ are birational for every $s\in S$.

\medskip
 {\it Warning \ref{bir.fam.defn}.1.}  If   each $\phi_s$ is an isomorphism, then 
$\phi$ is  an isomorphism over a dense, open subset of $S$ by (\ref{birm.mps.say}.2), 
 but  it can happen that $\phi$  is not an isomorphism, see (\ref{small.mod.rem}.2) for such an example.
\end{defn}

\begin{say}[Restriction of birational maps] \label{birm.mps.say}
  Let   $\phi:X\map Y$ be a birational map, and
 $W_X\subset X$ and $W_Y\subset Y$  closed subsets such that
  $\phi$ restricts to a birational map  $\phi_W:W_X\map W_Y$.
  Then $\dom(\phi)\cap W_X\subset \dom(\phi_W)$, but they need not be equal.
That is, it can happen that $\phi_W$ is an isomorphism at a point $x\in W_X$, but
$\phi$ is not  defined at $x$.

It is  useful to know that this does not happen if
$W_X$ is in general position.

 Let $S$ be a normal scheme and   $\pi:X\to S$ an irreducible, normal scheme, proper over $S$, with geometrically irreducible generic fiber.
There is a dense, open $S^\circ\subset S$ over which the fibers are
geometrically irreducible. Set $S_1:=S\setminus S^\circ$.
For $i\geq 2$ let $S_i\subset S$ be the closed subset over which the fiber dimension is $\geq i-1+\dim X-\dim S$.
By a dimension count we see the following.

\medskip
{\it Claim \ref{birm.mps.say}.1.}   
Let $W\subset S$ be an irreducible,  local complete intersection of codimension $m$.
If $\codim(W\cap S_i, S_i)\geq m$ for every $i$, then
$\pi^{-1}(W)$ is irreducible. \qed
\medskip

 Next let $\pi_i:X_i\to S$ be irreducible, normal schemes, proper over $S$, and $\phi:X_1\map X_2$ a birational map with closed graph
  $\Gamma(\phi)$. Applying (\ref{birm.mps.say}.1) to the $X_i\to S$ and to
  $\Gamma(\phi)\to S$, we get the following.
  
  \medskip
      {\it Claim \ref{birm.mps.say}.2.}
      There are closed subsets $S_i\subset S$ with the following property.
      
      Let $Z\subset S$ be an irreducible,  local complete intersection of codimension $m$ such that 
      $\codim(Z\cap S_i, S_i)\geq m$ for every $i$. Then
      $\phi$ restricts to a birational map
      $\phi_Z: \red\pi_1^{-1}(Z)\map  \red\pi_2^{-1}(Z)$, and
      $\dom(\phi)\cap \pi_1^{-1}(Z)= \dom(\phi_Z)$. \qed
      \end{say}

We are especially interested in cases when $\dom(\phi)$ is large.

\begin{defn}[Small  modification] \label{small.mod.defn}
An open subset $U\subset X$ is called {\it large} if
$X\setminus U$  has
codimension $\geq 2$ in $X$.
A birational map  $\phi: X_1\map X_2$ is a
{\it small  modification} if  $\dom(\phi)$ is large in both $X_1$ and $X_2$.

The composite of small  modifications is a small  modification.

A family of birational maps $\phi: X_1\map X_2$  over $S$ is a
 {\it family of small  modifications} if
 $\phi_s: (X_1)_s\map (X_2)_s$ is a small  modification for every $s\in S$.

\end{defn}

\begin{rem}\label{small.mod.rem}
  Given $X_i\to S$,  let  $\phi:X_1\map X_2$ be a small  modification.
  By definition there are  large open subsets   $U_i\subset X_i$  such that $\phi$ restricts to an
isomorphism  $\phi_U: U_1\cong U_2$. 
There is a dense, open subset $S^\circ\subset S$ such that
$(U_i)_s\subset (X_i)_s$ are large open subsets  for $s\in S^\circ$.

Thus $\phi$ is a  family of small  modifications over $S^\circ$, but not always over $S$. The following example illustrates some of the possibilities.

  Start with the deformation of the quadric cone
  $$
  X:=(x_1^2-x_2^2+x_3^2=t^2)\subset \p^3_x\times \a^1_t.
  $$
  It has 2 small simultaneous resolutions. We get $\pi_1:X_1\to X$ by blowing up
  $(x_1-x_2=x_3-t=0)$, and $\pi_2:X_2\to X$ by blowing up $(x_1+x_2=x_3-t=0)$. 
  \begin{enumerate}
  \item The $\pi_i:X_i\to X$ are small  modifications, but not
    families of small  modifications over $\a^1$.
    \item 
  The composite  $\phi:=\pi_2^{-1}\circ \pi_1: X_1\map X_2$ is a family of small  modifications,
  where each $\phi_t$ is an isomorphism. However $\phi$ is not an isomorphism.
\end{enumerate}
\end{rem}

Recall that a divisor $D\subset X$  is {\it $\pi$-exceptional}
for $\pi:X\to S$ if  $\pi(D)\subset S$ has codimension $\geq 2$.

\begin{lem} Let $S$ be a normal scheme,  $\pi_i:X_i\to S$ irreducible, normal schemes, proper over $X$, and $\phi:X_1\map X_2$ a birational map.
  Assume that there are no $\pi_i$-exceptional divisors.
  Then $\phi$ is a small modification iff $\phi_s:(X_1)_s\map (X_2)_s$ is
  \begin{enumerate}
  \item   a small modification for the generic point $s\in S$, and
  \item    birational  for each codimension 1  point $s\in S$. \qed
  \end{enumerate}
  \end{lem}

The following isomorphism criterion is proved in \cite{ma-mu};
see \cite[11.39]{k-modbook} for the relative version.

\begin{thm}[Matsusaka-Mumford criterion]\label{m-m.lem}  Let $\pi_i:X_i\to S$ be normal schemes, projective over $S$, and  $\phi:X_1\map X_2$  a small  modification.
  Let  $H_1$ be a $\pi_1$-ample $\q$-Cartier divisor. Assume that
  $H_2:=\phi_*(H_1)$  is   $\q$-Cartier and $\pi_2$-ample.

  Then $\phi$ is an isomorphism.\qed
  \end{thm}

Next we show that small, projective, normal  modifications  are rigid.
It is likely that the same holds for  proper  modifications, but
the proof uses projectivity.

\begin{thm}\label{small.mod.prop}
  Let $X, S$ be  geometrically irreducible and  geometrically normal $k$-varieties,    $Y\to S$  a  projective morphism with  geometrically normal generic fiber, and
  $\phi: X\times S\map Y$ a  small  modification.
  Assume that
  \begin{enumerate}
  \item either $\phi$ is a  family of birational maps,
  \item or $S$ has rational singularities.
    \end{enumerate}
  Then there is a unique, small, projective  modification  $\sigma:X\map X'$ such that the composite
  $$
  \phi\circ (\sigma, 1_S)^{-1}:  X'\times S\map  X\times S\map Y
  \qtq{is an isomorphism.}
  $$
\end{thm}

Example~\ref{small.mod.prop.2.exmp} shows that assumption (\ref{small.mod.prop}.2) is necessary.
\medskip

Proof. Assume first that the base field is algebraically closed.
As we noted in (\ref{birm.mps.say}.2), 
there is a dense, open subset $S^\circ\subset S$ such that
$Y\to S$  is flat with geometrically normal fibers over $S^\circ$,
and $\dom(\phi^{-1})\cap Y_s=\dom(\phi_s^{-1})$ for $s\in S^\circ$, and  $\phi$ is a  family of small  modifications over $S^\circ$.

 Pick $s\in S^\circ$ and  take
 $X':=Y_s$. The 
composite map
  $$
  \psi:=\phi\circ (\phi_s, 1_S)^{-1} :  X'\times S\map  X\times S\map Y
  $$
  is then a small  modification such that
  $\dom(\psi)\cap X'\times \{s\}\subset X'\times \{s\}\cong Y_s$ is large, and  
$\psi_s:  Y_s=X'\times \{s\}\map Y_s$ extends to the identity map.

  Let $H'_s\subset X'$ be a normal, ample divisor (\ref{div.norm.say}.1). 
  Set $H':=H'_s\times S$ and $H:=\psi_* (H')$  its
  birational transform.  Note that $H_s=H'_s$.
  Next  (\ref{div.norm.say.3})  says that $H$ is Cartier along $Y_s$.
  Since $H|_{X_s}=H'_s$ is ample,  $H$ is relatively ample  over $S^\circ$,  after shrinking $S^\circ$ as needed. Then  (\ref{m-m.lem})  shows that $\psi$ is an isomorphism over $S^\circ$.

Thus  $\psi: X'\times S\map Y$ is a small  modification, which is an
isomorphism over $S^\circ$. Let $A$ be an ample divisor on $Y$ and
$A':=\psi^{-1}_*(A)$ its birational transform.

Since $\psi$ is an
isomorphism over $S^\circ$, $A'$ is  Cartier and ample over  $S^\circ$.

If  $\phi$ is a  family of birational maps, then
$supp A'$ does not contain any fiber $X'_s$, so 
$A'$ is 
a  family of Weil divisors. Thus 
  $A'$ is  Cartier and ample over  $S$ by (\ref{div.norm.say}.2). 
We can now use  (\ref{m-m.lem})  again to get that  $\psi$ is an isomorphism.

Note that $A'$ gives a morphism   $S^\circ\to \pics(X')$, which extends to a
morphism   $S\to \pics(X')$ if $S$ has rational singularities by
(\ref{map.to.av.say}).
Then $A'$ is Cartier and ample over  $S$.
As before,  $\psi$ is an isomorphism by  (\ref{m-m.lem}). 

For an arbitrary  base field $k$, we can choose the above $s\in S(\bar k)$ to be Galois over $k$. Then $Y_{\bar k}\cong Y_s\times S_{\bar k}$, thus the conjugate fibers of $Y_{\sigma(s)}$ are isomorphic to $Y_s$. 
Therefore $Y_s$ is defined over $k$.
\qed

\begin{rem} Note that there may be several  ways of writing
  $Y$ as $X'\times S$. 
  For example, if $G$ is an algebraic group and $X$ is a nontrivial $G$-torsor,
  then  $G\times X\cong X\times X$.

  The claim of (\ref{small.mod.prop}) is that the rational product structure
  $\phi$ can be promoted to a regular  product structure in a unique way.
\end{rem}

\begin{say}[Proof of Theorem~\ref{sm.mod.gr.isom.cor}]\label{sm.mod.gr.isom.cor.say}
Take 2 copies  $Y_i:=X\times G$ and consider the
rational map
$$
\phi: Y_1\map Y_2\qtq{given by}  (x, g)\mapsto \bigl(\mu(g,x), g\bigr).
$$
Let $\Gamma\subset Y_1\times_G Y_2$ be its closed graph.
By (\ref{birm.mps.say}.2) there is a dense, open subset $U\subset G$ such that, for $g\in U$,  the fiber
$\Gamma_g$ equals the closed graph of $\mu_g:X\map X$.
Thus $\phi$ is a family of small  modifications over $U$ by (\ref{birm.mps.say}.2).

By (\ref{small.mod.prop})
   there is a small, projective  modification  $\sigma:X\map X'$ such that the composite
  $$
  \psi: (\sigma, 1_U)\circ\phi:  X\times U\map  X\times U \map X'\times U
  \qtq{is an isomorphism.}
  $$
  Thus,  we get that the composite
  $$
  \mu_h^{-1}\circ \mu_g= \psi_h^{-1}\circ \psi_g:  X\map X
  $$
  is an isomorphism for $g, h\in U$.
  Since $G$ is connected,  $(h,g)\mapsto h^{-1}g$
is a surjective map  $U\times U\to G$.
  Hence
$\mu_g$ is an isomorphism for  all
  $g\in G$.

  As we noted in  (\ref{bir.fam.defn}.1), this implies that
  $\phi$  is an isomorphism over a dense, open subset of $G$.
  The $G$-action now shows that it is an isomorphism  everywhere.
  \qed
\end{say}

\begin{say}[Divisors on normal varieties]\label{div.norm.say}
  We used two results about normal, projective varieties $X$.

  (\ref{div.norm.say}.1) Let $|H|$ be a very ample linear system. Then a general $H\in |H|$ is normal; see \cite[IV.12.1.6]{ega}.

  (\ref{div.norm.say}.2) $\pico(X)$ is proper. Thus if  $D\subset X\times S$ is a  family of Weil divisors that is Cartier along the generic  fiber, then it is Cartier.
  Over $\c$ this goes back to Picard and Severi, the general case is in \cite{mats-pic}; see \cite[Sec.11.3]{klos} for a detailed discussion.

 \end{say}

\begin{say}[Lifting normal, Cartier divisors]\label{div.norm.say.3}
  Let $Y$ be a normal scheme,  $X\subset Y$ a Cartier divisor and
  $D$ an effective Weil divisor on $Y$. Assume that the
  {\it divisorial restriction}
  $D|_X^{\rm div}$ is normal and Cartier. Then $D$ is Cartier along $X$.

  Recall that the divisorial restriction is obtained by first taking the scheme-theoretic restriction  $\o_D|_X$, and then the quotient by its torsion subsheaf; see \cite[p.106]{k-modbook} for details.

  Indeed, since  $D|_X^{\rm div}$ is normal, there is a closed $Z\subset X$ of codimension $\geq 3$ such that $D|_X^{\rm div}$  is regular on $X\setminus Z$.
  Since $D|_X^{\rm div}$ is Cartier, $X$ and hence $Y$ are regular along
  $(D\cap X)\setminus Z$.  Thus $D$ is  Cartier along $X\setminus Z$.

  Next note that  $\depth_Z\bigl(D|_X^{\rm div}\bigr)\geq 2$ since  $D|_X^{\rm div}$ is normal, hence
  $\depth_Z(X)\geq 3$ since  $D|_X^{\rm div}$ is Cartier. 
  Thus  $D$ is Cartier along  $Z$ by  \cite[XI.2.2]{sga2}; see also
  \cite[10.73]{k-modbook}.
  \end{say}

\begin{say}[Rational maps to Abelian varieties] \label{map.to.av.say}
  Let $X$ be a normal variety, $A$ an Abelian variety and
  $\phi:X\map A$ a rational map. We would like to know that it is necessarily a morphism.

  If $X$ is smooth, this is  \cite[II.15.Prop.1]{MR29522}. For the general case, let $Y\to X\times A$ be the normalization of the closed graph of $\phi$.  Let $L$ be the pull-back of the universal  bundle from $A\times \pico(A)$ to
  $Y\times \pico(A)$.

  If $p:Y\to X$ is not an isomorphism, then there is a curve $C\subset Y$
  whose second projection to $A$ is not constant. Thus the restriction of $L$ to $C \times \pico(A)$ is nonconstant.

  On the other hand, the restriction of $L$ to $Y\times [\o_A]$ is $\o_Y$.   Since $R^1p_*\o_Y=0$, the constant 1 secion of $\o_Y$ lifts to a
  section of $L$ in a neighborhood of $[\o_A]\in \pico(A)$. This is a contradiction. \qed

  \end{say}
  
We also need a  scheme theoretic version of (\ref{small.mod.prop}).
The proof is essentially the same, but   seems to need more assumptions and
different references.
Also, it is not clear to me how to descend projectivity  to arbitrary base schemes from the
strictly Henselian case, so we state  only the latter.

 \begin{thm}\label{small.mod.prop.gen}
   Let $S$ be a normal, local, strictly Henselian scheme, $X$ a normal $S$-scheme,
   and  $W\to S$ a smooth morphism with geometrically connected fibers.
     Let  $Y\to W$  be a  projective morphism and
  $\phi: X\times_S W\map Y$ a  small   modification.
  Assume that
  \begin{enumerate}
  \item $Y_{\tau}:=S\times_{\tau, W}Y$ is normal for every section $\tau:S\to W$, and
  \item the induced map  $\phi_{\tau}: X\map Y_{\tau}$ is a small modification.
    \end{enumerate}
  Then there is a small, projective   modification  $\sigma:X\map X'$ such that the composite
  $$
  \phi\circ (\sigma, 1_W)^{-1}:  X'\times_S W\map  X\times_S W\map Y
  \qtq{is an isomorphism.}
  $$
\end{thm}

 Proof.  Note that section of $W\to S$ are local complete intersections. Thus 
 by (\ref{birm.mps.say}.2)  there is section $\tau:S\to W$ such that
 $\dom(\phi^{-1})\cap Y_{\tau}=\dom(\phi_{\tau}^{-1})$. 
 Set $X':=Y_{\tau}$ and consider the
composite map
  $$
  \psi:=\phi\circ (\phi_{\tau}, 1_W)^{-1} :  X'\times_S W\map  X\times_S W\map Y.
  $$
  Then $\dom(\psi)$ contains a large open subset of $X'\times_S\tau(S)$, and  
  $$\psi_{\tau}: X'\times_S\tau(S)\cong X'= Y_{\tau}\map Y_{\tau}$$ extends to the identity map.

  Let $H'_{\tau}\subset X'$ be a normal, relatively ample divisor (\ref{div.norm.say}.1). 
  Set $H':=H'_{\tau}\times_S W$ and $H:=\psi_* (H')$  its
  birational transform on $Y$.  Note that $H_{\tau}=H'_{\tau}$.
  As before,   (\ref{div.norm.say.3})  implies that
  $H$ is Cartier along $Y_{\tau}$.

  Since $H|_{X_{\tau}}=H'_{\tau}$ is  relatively  ample,  $H$ is relatively ample  over
  a suitable open subset  $\tau(S)\subset W^\circ\subset W$. Then  (\ref{m-m.lem})  shows that $\psi$ is an isomorphism over $W^\circ$, and $W^\circ$ dominates $S$ since it contains a section.

Thus  $\psi: X'\times_S W\map Y$ is a small  modification, which is an
isomorphism over $W^\circ$. Let $A$ be an ample divisor on $Y$ and
$A':=\psi^{-1}_*(A)$ its birational transform.

Since $\psi$ is an
isomorphism over $W^\circ$, $A'$ is  Cartier and ample over  $W^\circ$.
Thus  $A'$ is  Cartier  by (\ref{ram-sam.thm}),  and ample over  $W$. 
We can now use  (\ref{m-m.lem})  again to get that  $\psi$ is an isomorphism.
 \qed

\begin{prop}\label{gp.sch.acts.prop}
   Let $S$ be a normal scheme, $X$ a normal, projective $S$-scheme,
   and  $G\to S$ a smooth group scheme with  irreducible fibers.
   Let  $\mu_S: G\times_S X\map X$ be a rational action.
   Assume that for every \'etale section  $\tau:S'\to G$  (that is,  $S'\to G\to S$ is  \'etale), the induced map
   $$
   \tau_X:X_{S'}=(S'\times_S X)\stackrel{(1_{S'},\tau, 1_X)}{\longrightarrow}
   S'\times_SG\times_S X\stackrel{(1_{S'},\mu_S)}{\longrightarrow}
   (S'\times_S X)=X_{S'}
   \eqno{(\ref{gp.sch.acts.prop}.1)}
   $$
  is a small modification. 
   Then $\mu_S: G\times_S X\map X$ is a regular action.
\end{prop}

Proof.  Being a morphism is an \'etale-local property, we may thus assume that
$S$ is local and strictly Henselian.

Take 2 copies  $Y_i:=X\times_S G$ and consider the
rational map
$$
\phi: Y_1\map Y_2\qtq{given by}  (x, g)\mapsto \bigl(\mu(g,x), g\bigr).
$$
By (\ref{gp.sch.acts.prop}.1) the  assumption   (\ref{small.mod.prop.gen}.2) is satisfied. 
Thus by (\ref{small.mod.prop.gen})
   there is a small, projective  modification  $\sigma:X\map X'$ such that the composite
  $$
  \psi: (\sigma, 1_G)\circ\phi:  X\times_SG\map  X\times_SG \map X'\times_SG
  \eqno{(\ref{gp.sch.acts.prop}.2)}
  $$
  is an isomorphism. Let $\tau_i:S\to G$ be general sections.
  Restricting (\ref{gp.sch.acts.prop}.2) to $\tau_i(S)$ and composing, 
    we get that 
  $$
  \mu_{\tau_1}^{-1}\circ \mu_{\tau_2}= \psi_{\tau_1}^{-1}\circ \psi_{\tau_2}:  X\map X
  $$
  is an isomorphism.
  Since $G$ has  connected fibers, these $\mu_{\tau_1}^{-1}\circ \mu_{\tau_2}$ generate it.   Thus, by (\ref{birm.mps.say}.2), $\mu_S$ is an isomorphism over an open dense set
  of $G$ that contains a section.  By the $G$-action, it is an
  isomorphism everywhere.  \qed

\begin{say}[Proof of Theorem~\ref{main.min.ner.thm.gen}]\label{main.min.ner.thm.gen.pf}
Let  $\tau:S'\to G$ be an  \'etale section. As in
(\ref{gp.sch.acts.prop}.1) we get a birational map 
$\tau_X:X_{S'}\map  X_{S'}$.
Also,  $X_{S'}$  has terminal singularities and $K_{X_{S'}}/S$ is reletively nef.
Thus $\tau_X$ is a small modification by (\ref{min.mod.say}.4).
Thus $\mu_S: G\times_S X\map X$ is regular by (\ref{gp.sch.acts.prop}). \qed
\end{say}
\medskip

The following can be viewed as  a generalization of (\ref{div.norm.say}.2); for proofs  see \cite[IV.21.14.1]{ega} or \cite[10.65]{k-modbook}.

\begin{thm}[Ramanujam-Samuel]\label{ram-sam.thm}  Let $Z$ be a normal scheme  and
  $g:X\to Z$ a smooth morphism. Let $D$ be a Weil divisor on $X$
that is Cartier on a dense open subset of every  fiber of $g$.
   Then $D$ is Cartier. \qed
\end{thm}

  \section{Equivariant completions}\label{eq.comp.sec}

  We prove Corollary~\ref{main.min.ner.thm.cor.2} in 2 steps.

\begin{say}[Proof of \ref{main.min.ner.thm.cor.2}.1]\label{eq.comp.step.1}
This part uses characteristic 0.

  Let $\pi:X_T\to T$ be a $\q$-factorial minimal model of $A_K$, and  $H_0$  a $\pi$-ample divisor on $X_T$.

   By  (\ref{min.mod.say}.4)
     $\ner(A_K)$ acts on the N\'eron-Severi group of $X_T$,  and $\nero(A_K)$ acts trivially.
    So the $\ner(A_K)$-orbit of $[H_0]\in \ns(X_T)$ is finite. Choose representatives $\{H_i:i\in I\}$ from each divisor class in the orbit, and set  $H:=\sum_{i\in I}H_i$. Let
    $X^{\rm c}_T\to T$ be the canonical model of
    $\bigl(X_T,  H\bigr)$. That is
    $$
   X^{\rm c}_T:=\proj_T \oplus_m \pi_*\o_{X_T}(mH),
   $$
   provided  this sheaf of algebras is finitely generated.
    A technical problem is that minimal model theory is usually formulated for
   effective  divisors, and it is not clear that the $H_i$ can be chosen effective.
 However,  finite generation can be checked after an \'etale base change by \cite[6.6]{km-book}, so
   it is enough to show that $X^{\rm c}_T$ exists when $T$ is strictly Henselian.

   In this case there are sections $\sigma_i:T\to \ner(A_K)$
   passing through each irreducible component of the central fiber.
   We can then take  $H_i:=\sigma_i(H_0)$ and choose 
   $0<\epsilon\ll 1$  such that
   $\bigl(X_T,  \epsilon H\bigr)$ is log terminal.  Then
   $X^{\rm c}_T $ is its canonical model by \cite[3.52]{km-book}, and  $X_T\map X^{\rm c}_T $ is a small modification.

   Then $H$ descends to a relatively  ample divisor  $H^{\rm c}$ on $X^{\rm c}_T$, and  its class  $[H^{\rm c}]\in \ns(X^{\rm c}_T)$ is
   $\ner(A_K)$-invariant. Thus the  $\ner(A_K)$-action is regular on $X^{\rm c}_T$ by (\ref{m-m.lem}).
  \end{say}

\begin{say}[Proof of \ref{main.min.ner.thm.cor.2}.1 $\Rightarrow$  \ref{main.min.ner.thm.cor.2}.2]\label{eq.comp.step.2}
The arguments  work over any 1-dimensional,  regular base scheme.

  We have $Y_0:=X^{\rm c}_T\to T$ with  a $\ner(A_K)$-action and a section $e_0:T\to Y_0$, giving
    $p_0:\ner(A_K)\to  Y_0$ that is an isomorphism on the generic fiber.

    If $ Y_i$ is aready defined with a $\ner(A_K)$-action and a section $e_i:T\to Y_i$, let $p_i:\ner(A_K)\to  Y_0$ be the corresponding morphism.

    Let $Z_i\subset Y_i$ be the closure of the image of the central fiber
    $\ner(A_K)_k$. If $Z_i$ is a divisor and $X_i$ is regular at the generic points of $Z_i$, then
    $p_i$ restricts to an open embedding  $\ner(A_K)_k\into Z_i$, and
    so $p_i:\ner(A_K)\to  Y_0$ is an open embedding. Thus the  normalization of 
    $Y_i$ is an equivariant completion of  $\ner(A_K)$.

    Otherwise let $Y_{i+1}:=B_{Z_i}Y_i$ be the blow-up of $Z_i$.
    Since $Z_i$ is $\ner(A_K)$-equivariant, the
    $\ner(A_K)$-action on $Y_i$ lifts to a $\ner(A_K)$-action on $Y_{i+1}$,
    and the section $e_i$ also lifts to $e_{i+1}:T\to Y_{i+1}$.

    By a lemma of Abhyankar and Zariski  (cf.\ \cite[2.45]{km-book}), the process eventually stops, and
    we get the required equivariant completion of  $\ner(A_K)$.

    \end{say}

\section{Semiabelian reduction}\label{semiab.sec}

The correspondence between the minimal and N\'eron models is especially strong if $A_K$ has  {\it semiabelian reduction,} that is, the central fiber of $\nero(A_K)$ is an extension of an Abelian variety by a torus.

One can work with {\it log terminal minimal models.} That is,
(\ref{min.mod.say}.2--3) hold, but $X_T$ is allowed to have log terminal singularities.
 Using  minimal model  theory, we get the following in  characteristic 0.

\begin{thm} \label{loc.st.char.thm}
  Let   $A_K$  be  an  Abelian variety,
  and  $\pi:X_T\to T$  a  log terminal minimal model  of $A_K$ with  central fiber  $X_k$.
  Let $A_k$ denote the identity component of the central fiber of $\ner(A_K)$.
   The following are equivalent.
 \begin{enumerate}
 \item  $A_k$   is a semiabelian variety.
   \item  $X_k$  is reduced, and  the pair $(X_T, X_k)$ is log canonical.
 \end{enumerate}
\end{thm}

If (\ref{loc.st.char.thm}.2) holds then  $\pi:X_T\to T$  is {\it locally stable} in the terminology of \cite{k-modbook}.
\medskip
We start the proof with some auxiliary lemmas.

\begin{lem}\label{red.comp.diff.lem}
  Let $\pi:Y_T\to T$ be a proper morphism such that $Y_T$ is normal and
$K_{Y_T}$ is $\q$-Cartier.  Let  $E\subset Y_k$ be an irreducible component
  of multiplicity 1.  Let $E^*\subset E$ be the open set where $Y_k$ is smooth, and $D\subset E$ its complement.
Let $\tau: \bar E \to E$ be  normalization and 
$\bar D\subset \bar E$ the divisorial part of the preimage of $D$.

Then there  is an effective $\q$-divisor $D^*$ suported on $\bar D$ such that
$$
K_{\bar E}+\bar D+D^* \simq \tau^*\bigl(K_{Y_T}+Y_k)\simq \tau^*K_{Y_T}.
\eqno{(\ref{red.comp.diff.lem}.1)}
$$
\end{lem}

Proof.   Write  $\Delta:=Y_k-E$, as a divisor on $Y_T$.
By \cite[4.2]{kk-singbook}, there is an effective divisor
$\diff_{\bar E} \Delta $  such that
$$
K_{\bar E}+\diff_{\bar E} \Delta  \simq \tau^*\bigl(K_{Y_T}+Y_k)\simq \tau^*K_{Y_T},
\eqno{(\ref{red.comp.diff.lem}.2)}
  $$
Let $D_i\subset D$ be an irreducible component, and $\bar D_i\subset \bar D$ its preimage.

If $Y_k$ is singular along $D_i$, then  $(Y_T, Y_k)$ is not klt along $D_i$, hence $\bar D_i$ appears in
$\diff_{\bar E} \Delta $ with coefficient $\geq 1$ by \cite[2.35.5--6]{kk-singbook}.   We can thus set $D^*:=\diff_{\bar E} \Delta-\bar D$.  \qed

\medskip

The following is a special case of \cite{MR3793359}.

\begin{lem} \label{eq.comp.cases.lem}
  Let $A$ be a connected, commutative algebraic group with unipotent radical $U\subset A$. Let 
  $\bar A\supset A$ be a normal, $A$-equivariant compactification, and $\bar D\subset \bar A$ the  divisorial part of  $\bar A\setminus A$.
    \begin{enumerate}
      \item If $U\neq \{0\}$, then $K_{\bar A}+\bar D$ is not pseudo-effective.
      \item If $U= \{0\}$, then $K_{\bar A}+\bar D\sim 0$ and  $(\bar A, \bar D)$ is log canonical. \qed
    \end{enumerate}
\end{lem}

\begin{say}[Proof of (\ref{loc.st.char.thm}.2) $\Rightarrow$ (\ref{loc.st.char.thm}.1)]\label{lc.to.can.cor}
   We check in  (\ref{logmin.lc.thm}) that there is a regular $\nero(A_K)$-action on $X_T$.

Apply (\ref{red.comp.diff.lem}) to  $Y_T:=X_T$ and any irreducible component $E\subset X_k$.
Note that  $E^*\subset E$ is the dense $A_k$-orbit by (\ref{sm=orbit.lem}). 

Then (\ref{red.comp.diff.lem}.2) gives that
$K_{\bar E}+\diff_{\bar E} \Delta\simq 0$. 
Since  $(X_T, X_k)$ is log canonical, 
$\diff_{\bar E}\Delta=\bar D$  by \cite[4.9]{kk-singbook}.
Thus $K_{\bar E}+\bar D\simq 0$, hence $A_k$ is semiabelian
by (\ref{eq.comp.cases.lem}.1). \qed
  \end{say}

\begin{say}[Proof of (\ref{loc.st.char.thm}.1) $\Rightarrow$ (\ref{loc.st.char.thm}.2)]\label{loc.st.char.thm.pf.2}

As in (\ref{eq.comp.step.2}),  let
$Y_T\to T$ be a   $\nero(A_K)$-equivariant resolution of $X_T$ that contains
$\ner(A_K)$ as an open subset.
Running a  minimal model program over $X_T$ gives
$Y_T=Y_T^1\map \cdots \map Y_T^m=:Z_T$.

Let $E\subset Y_k$ be any divisor corresponding to an irreducible component of
 $\ner(A_K)_k$.
We claim that $E$ is not contracted by
$Y_T\map Y_T^m$.

Assume by induction that $E$ appears on $Y_T^i$. 
If $A_k$ is semiabelian then  $K_{\bar E}+\bar D\sim 0$ by (\ref{eq.comp.cases.lem}.2),
hence
$
\tau^*K_{Y_T^i}\simq K_{\bar E}+\bar D+D^*\simq D^*
$
is effective. Thus $Y_T^i\to Y_T^{i+1}$  does not contract $E$.

Therefore  $Z_k$ contains at least 1 divisor $E$ with multiplicity 1. 
Here  $K_{Z_T}\simq 0$, hence 
$K_{\bar E}+\diff_{\bar E} \Delta \simq 0$, and so
$\diff_{\bar E}\Delta=\bar D$. 
Inversion of adjunction \cite[4.9]{kk-singbook}
implies that $(Z_T, Z_k)$ is log canonical along $E$.
In particular, every irreducible component of $Z_k$ meeting $E$ also  has
multiplicity 1. Repeating the argument for the other
irreducible components we get that $Z_k$ is reduced and 
$(Z_T, Z_k)$ is log canonical. Then the same holds for 
 $(X_T, X_k)$. 
\qed
\end{say}

\medskip

\section{Log canonical models}\label{logmin.sec}

In this section 
$S$ denotes a normal, integral  base scheme with function field $K$.
The proof   uses several results of the minimal model program, so we need to assume that $S$ is excellent, has a dualizing complex, and is  over a field of characteristic 0.
With these restrictions, 
regularity of actions  also holds for almost all pairs  over higher dimensional bases.

 \begin{thm}\label{gp.sch.acts.mm.prop}
   Let 
   $p:(X_S,\Delta_S)\to S$ be a dominant, projective morphism,
   and  $G\to S$ a smooth group scheme with  irreducible fibers.
   Let $\mu_K:G_K\times X_K\to X_K$ be a regular action such that
   $\Delta_K$ is $G_K$-invariant.
   Assume that
   \begin{enumerate}
   \item $(X_S,\Delta_S)$ is log canonical and $(X_K,\Delta_K)$ is terminal,
      \item  $K_{X_S}+\Delta_S$ is $p$-semiample, and
   \item  $  p^{-1}\bigl(p(Z)\bigr)\subset \supp \Delta_S$
   for every log canonical center $Z$ of $ (X_S,\Delta_S)$.
   \end{enumerate}
   Then the induced rational action
     $\mu_S: G\times_S X\map X$ is regular.
\end{thm}

 {\it Remarks \ref{gp.sch.acts.mm.prop}.4.}
See \cite[4.28]{kk-singbook} for {\it log canonical centers.}
 
Assumptions (\ref{gp.sch.acts.mm.prop}.1--2) generalize the notion of a minimal model to pairs (allowing worse singularities than usual).

The abundance conjecture says that if $K_{X_S}+\Delta_S$ is $p$-nef, then it is also $p$-semiample,
This is known if $(X_S,\Delta_S)$ is klt and $K_{X_K}+\Delta_K\simq 0$ \cite{fujino-ssmmp, lai-2009}. 

Let $E$ be a smooth, elliptic curve and $X$ the blow-up of a point
$(e,0)\in E\times \a^1$. Let $\Delta:=E'$ be the birational transform of $E_0$ 
with coefficient 1.  Then $K_X+\Delta$ is numerically trivial but the
$E$-action does not lift to $X$. Note that  $\Delta$ has a component
 with coefficient 1 (namely $E'$) and also  a component with coefficient 0 (the exceptional $\p^1$).
Assumption (\ref{gp.sch.acts.mm.prop}.3) is there to exclude such extreme cases.

\medskip

Proof. Let $g: X'_S\to X_S$ be a $\q$-factorial,
minimal modification. That is, 
 $X'_S$ has terminal singularities, and $K_{X'_S}$ is $g$-nef.
Thus there is a unique effective $\q$-divisor $\Delta'_S$
such that $K_{X'_S}+\Delta'_S$ is numerically $g$-trivial and
$g_*\Delta'_S=\Delta_S$.
These imply that  $(X'_S,\Delta'_S)$ is log canonical and
$\supp\Delta'_S\supset g^{-1}(\supp\Delta_S)$ by \cite[3.39]{km-book}.

By (\ref{gp.sch.acts.mm.prop}.2) there is a  factorization $p:X\to S'\to S$
such that $S'\to S$ is birational and $K_{X_S}+\Delta_S$ is numerically trivial on $X_S\to S'$.
MMP for $X'_S\to S'$ gives a minimal model  $X^{\rm m}\to S'$.
Let $\tau:X'_S\map X^{\rm m}$ be the corresponding birational map.

At each step of the MMP the image of $K_{X'_S}+\Delta'_S$ is numerically trivial, and
we end with  a log canonical pair $(X^{\rm m}, \Delta^{\rm m})$ such that
$K_{X^{\rm m}}+ \Delta^{\rm m}$ is also numerically trivial over $S'$.
  For any irreducible divisor  $E\subset  X'_S$, 
  $$
  a(E, X^{\rm m}, \Delta^{\rm m})=a(E, X'_S, \Delta'_S)=a(E, X_S, \Delta_S)
  $$
  is the negative of its
  coefficient in $\Delta'_S$.
  Thus, using the notation in (\ref{totald.defn}),
  $$
  a(E, X^{\rm m}, \Delta^{\rm m})\leq 0\leq 1+\totaldiscrep\bigl(\tau(E),X^{\rm m}, \Delta^{\rm m}\bigr),
  $$
  and equality holds iff
  $$
  a(E, X^{\rm m}, \Delta^{\rm m})=0\qtq{and}
  \totaldiscrep\bigl(\tau(E),X^{\rm m}, \Delta^{\rm m}\bigr)=-1.
  $$
  The first says that $E$ is not contained in the support of $\Delta'_S$, the second that $g(E)$ is contained in a log canonical center of $(X_S, \Delta_S)$.
  This is excluded by  (\ref{gp.sch.acts.mm.prop}.3).

  Note that  $ X^{\rm m}$ is a minimal model over $S$, so
  the induced $G$-action on $ X^{\rm m}$ is regular by
  Theorem~\ref{main.min.ner.thm.gen}. It then lifts to a regular $G$-action on $ X'_S$ by (\ref{aut.lift.lem}),
  which pushes forward to a regular $G$-action on $ X_S$ by    (\ref{act.pushfwd.cor}).\qed

  \begin{lem}\label{act.pushfwd.cor}
    Let $S$ be a normal scheme, $X, Y$  normal, projective $S$-schemes,
   and  $G\to S$ a smooth group scheme with  irreducible fibers.
   Let  $\mu_X: G\times_S X\to X$ be a regular action
   and $g:X\map Y$ a birarational map. Assume that
   \begin{enumerate}
   \item   $g^{-1}$ has no exceptional divisors, and
   \item every $g$-exceptional divisor is $\mu_X$-invariant.
   \end{enumerate}
   Then $\mu_X$ pushes forward to a  regular action
   $\mu_Y: G\times_S Y\to Y$.
   \end{lem}

Proof. We have a rational action
   $\mu_Y: G\times_S Y\map Y$. For every \'etale section  $\tau:S'\to G$  the induced map
   $\tau_Y:=g\circ \tau_X\circ g^{-1}$
is a small modification by the assumptions. Thus $\mu_Y$ is a
regular action by (\ref{gp.sch.acts.prop}). \qed

\begin{defn}\label{totald.defn} Let $(X, \Delta) $ be a pair and  $W\subset X$ an irreducible subset. Let $\totaldiscrep(W, X, \Delta)$ denote the infimum of
  the discrepancies $a(E, X, \Delta)$, where $E$ runs through all divisors
  for which $W\subset \cent_XE$.
\end{defn}

\begin{lem}\label{aut.lift.lem}
Let $S$ be a normal scheme,
 $p:(X, \Delta)\to S$  a projective morphism, and
  $G\to S$ a smooth group scheme with connected fibers acting on $X$.

Let $g:Y\to X$ be a projective, birational  morphism such  that $Y$ is normal,
and
$                                                                              
a(E, X, \Delta)<1+\totaldiscrep(g(E), X, \Delta)                                      
$
for every $g$-exceptional divisor $E$.

Then the $G$-action on $X$ lifts to a  $G$-action on $Y$.
\end{lem}

Proof. Let $h:X'\to X$ be a $G$-equivariant log resolution of $(X, \Delta)$.
There is a unique snc dvisor $\Delta'$ such that
$K_{X'}+\Delta'\simq h^*(K_{X}+\Delta)$ and
$h_*\Delta'=\Delta$. Then
$a(E, X', \Delta')=a(E, X, \Delta)$ by \cite[2.30]{km-book}.  Since $a(E, X, \Delta)<1+\totaldiscrep(g(E), X, \Delta) $,
 the center of $E$ on $X'$ is a stratum of  $(X', \Delta')$  by
 \cite[2.31]{km-book}.
All strata are $G$-equivariant. Thus, after further $G$-equivariant blow-ups
$X''\to X'$, we may assume that  every $g$ exceptional divisor $E$ is a divisor
 on $X''$  by \cite[2.45]{km-book}. Then the $G$-action on $X''$ pushes forward
 to a  $G$-action on $Y$ by (\ref{act.pushfwd.cor}). \qed

 \medskip
Specializing (\ref{gp.sch.acts.mm.prop}) to the case when $S=T$ is the spectrum of a DVR and $X_K$ is an Abelian torsor gives the following. 

\begin{cor}\label{logmin.lc.thm}
   Let  $X_K$  be a torsor under an  Abelian variety $A_K$,
   and  $\pi:(X_T, \Delta)\to T$  a projective morphism.
   Assume that
   \begin{enumerate}
     \item either $(X_T, \Delta)$ is log canonical and 
       $\supp\Delta=X_k$,
     \item or $(X_T, \Delta)$ is klt and  $\supp\Delta\subset X_k$.
       \end{enumerate}
     Then the $A_K$-action on $X_K$ extends to a regular $\nero(A_K)$-action on $X_T$. \qed
\end{cor}

\begin{cor}\label{logmin.lc.thm.cor}   Continuing with the notation of (\ref{logmin.lc.thm}), assume in addition that $\Delta=\red (X_k)$,
  $(X_T, \Delta)$  is divisorial log terminal, and
  $X_k$ has an irreducible  component of  multiplicity 1.
  Then  
\begin{enumerate}
\item $\ner(X_K)$ is isomorphic to the $\pi$-smooth locus of $ X_T$, and
\item $\ner(A_K)_k$  is semiabelian.
  \end{enumerate}
\end{cor}

Proof.  We follow the proof of (\ref{gp.sch.acts.mm.prop}).
For $g:X'_T\to X_T$ we choose  a
{\it thrifty modification} as in \cite[2.79]{kk-singbook}.
That is,  all $g$-exceptional divisors have coefficient $<1$ in $\Delta'$.  Note that
$$
K_{X'_T}\simq -\Delta'\simq X'_k-\Delta'.
$$
If $E\subset X'_k$ is an irreducible component that dominates a
reduced component of $X_k$, then $E$ appears in $ X'_k-\Delta'$
with coefficient $0$, all other irreducible components appear with
positve coefficient. Therefore the  $K_{X'_T}$-MMP
contracts precisely the support of $ X'_k-\Delta'$.

By (\ref{G.circ.acts.thm.cor}.5) the  smooth locus of $X^{\rm m}_T\to T$ is isomorphic to $\ner(X_K)$,
hence the same holds for $X_T\to T$ by  (\ref{sm=orbit.lem}).
The second claim follow from  (\ref{loc.st.char.thm}). \qed

  \begin{exmps}\label{log.mm.exmps}
    These examples show that there are many log minimal models, and
    the $\supp\Delta=Y_k$  assumption in (\ref{logmin.lc.thm})  is needed.

\smallskip
    (\ref{log.mm.exmps}.1) Let $E$ be an elliptic curve.
    Blowing up a point in the central fiber of
    $E\times \a^1\to \a^1$, and then contracting the birational transform of $E$ gives a log minimal model  $\pi:S\to \a^1$.
    A pair  $(S, \Delta)$ is log canonical only for $\Delta=0$.
    There is no regular $\nero(A_K)$-action on $S$.
    \smallskip

    (\ref{log.mm.exmps}.2) Start with
    $
    \bigl(xyz=t(x^3+y^3+z^3)\bigr)\subset \p^2_{xyz}\times \a^1_t.
    $
    All log minimal models with $\Delta=X_k$ are constructed as follows.

    Repeatedly blow up singular points of the central fiber, then contract any
    subset of the central fiber (except we can not contract the whole cenral fiber).

    \smallskip

    (\ref{log.mm.exmps}.3) If in (\ref{log.mm.exmps}.2) we contract all but 1 or 2 of the irreducible components of the central fiber, we get a
    log minimal model that does not have a $\ner(A_K)$-action.
    
 \end{exmps}


  \newcommand{\etalchar}[1]{$^{#1}$}
\def\cprime{$'$} \def\cprime{$'$} \def\cprime{$'$} \def\cprime{$'$}
  \def\cprime{$'$} \def\dbar{\leavevmode\hbox to 0pt{\hskip.2ex
  \accent"16\hss}d} \def\cprime{$'$} \def\cprime{$'$}
  \def\polhk#1{\setbox0=\hbox{#1}{\ooalign{\hidewidth
  \lower1.5ex\hbox{`}\hidewidth\crcr\unhbox0}}} \def\cprime{$'$}
  \def\cprime{$'$} \def\cprime{$'$} \def\cprime{$'$}
  \def\polhk#1{\setbox0=\hbox{#1}{\ooalign{\hidewidth
  \lower1.5ex\hbox{`}\hidewidth\crcr\unhbox0}}} \def\cdprime{$''$}
  \def\cprime{$'$} \def\cprime{$'$} \def\cprime{$'$} \def\cprime{$'$}
\providecommand{\bysame}{\leavevmode\hbox to3em{\hrulefill}\thinspace}
\providecommand{\MR}{\relax\ifhmode\unskip\space\fi MR }
\providecommand{\MRhref}[2]{%
  \href{http://www.ams.org/mathscinet-getitem?mr=#1}{#2}
}
\providecommand{\href}[2]{#2}

  \bigskip

  Princeton University, Princeton NJ 08544-1000, \

  \email{kollar@math.princeton.edu}

\end{document}